\documentclass[12pt,a4letter]{amsart}
\usepackage[utf8]{inputenc}
\usepackage[english]{babel}

\usepackage{amsmath,amsxtra,amssymb,latexsym, amscd,amsthm}
\usepackage[mathscr]{eucal}

\newtheorem {theorem}{Theorem}[section]

\newtheorem {lemma}{Lemma}[section]
\newtheorem {example}{Example}[section]
\newtheorem {definition}{Definition}[section]
\newtheorem {remark}{Remark}[section]
\newtheorem {proposition}{Proposition}[section]

\newcommand{\R}{\mathbb{R}}

\def\EES{{\accent"5E e}\kern-.5em\raise.8ex\hbox{\char'23 }}
\def\ow{o\kern-.42em\raise.82ex\hbox{
   \vrule width .12em height .0ex depth .075ex \kern-0.16em \char'56}\kern-.07em}
\def\OW{o\kern-.460em\raise1.36ex\hbox{
\vrule width .13em height .0ex depth .075ex \kern-0.16em
\char'56}\kern-.07em}
\def\DD{D\kern-.7em\raise0.4ex\hbox{\char '55}\kern.33em}

\title{\L ojasiewicz-type inequalities and global error bounds for nonsmooth definable functions in o-minimal structures}

\author{Ho\`ang Phi D\~ung} 
\address{Department of Scientific Fundamentals,
Posts and Telecommunications Institute of Technology,
Office A2, Fl.10, Km10 Nguyen Trai Rd., Ha Dong District, Hanoi, Vietnam}
\email{dunghp@ptit.edu.vn}



\begin{document}
\maketitle

\begin{abstract}
In this paper, we give some \L ojasiewicz-type inequalities and a nonsmooth slope inequality on non-compact domains for continuous definable functions in an o-minimal structure. We also give a necessary and sufficient condition for which global error bound exists. Moreover, we point out the relationship between the Palais-Smale condition and this global error bound. 
\end{abstract}
\pagestyle{plain}

\section{Introduction}
Let $f : \mathbb{R}^n \to \mathbb{R}$ be a real analytic function with $f(0) = 0$. Let $V := \{x \in \mathbb{R}^n | f(x) = 0\}$ and $K$ be a compact subset in $\mathbb{R}^n$. Then the (classical) \L ojasiewicz inequality (see \cite{L1, L2}) asserts that:
\begin{itemize}
\item There exist $c > 0, \alpha > 0$ such that \begin{equation}\label{Eqn01}|f(x)| \ge cd(x, V)^\alpha\quad \text{for}\ x \in K.\end{equation}
\end{itemize}
\noindent Let $f : \mathbb{R}^n \to \mathbb{R}$ be a real analytic function with $f(0) = 0$ and $\nabla f(0) = 0$. The \L ojasiewicz gradient inequality (see \cite{L1, L2}) asserts that:
\begin{itemize}
\item There exist $C > 0, \rho \in [0, 1)$ and a neighbourhood $U$ of $0$ such that \begin{equation}\label{Eqn02} \|\nabla f(x)\| \ge C|f(x)|^\rho\quad \text{for}\ x \in U.\end{equation}
\end{itemize}
As a consequence, in (\ref{Eqn01}), the order of zero of an analytic function is finite, and if $f (x)$ is close to $0$ then $x$ is close to the zero set of $f$. However, if $K$ is not compact, the latter is not always true and the inequality (\ref{Eqn01}) does not always hold (see \cite[Remark 3.5]{DHN}). Similarly, in (\ref{Eqn02}), the order of gradient's zero of an analytic function is smaller than the order of its zero. But if ${U}$ is not a bounded set, (\ref{Eqn02}) does not always hold (see Example \ref{Ex2}).

With the \L ojasiewicz inequality (\ref{Eqn01}), in the case $K = \mathbb{R}^n$, H\"{o}rmander (see \cite{Hor}) substituted the left-hand side by one quantity greater than $|f(x)|$ and he got the following fact $$ \exists c, \alpha, \beta > 0\ \text{such that}\ |f(x)|(1 + |x|^\beta) \ge d(x, V)^\alpha, \forall x \in K.  $$ 

Recently, by replacing $V$ by a large real algebraic set, the authors in \cite{HN} and the authors in \cite{DHN} gave some versions of \L ojasiewicz inequalities in some non-compact cases. Moreover, some necessary and sufficicent conditions for which the \L ojasiewicz inequality and the global \L ojasiewicz inequality exists in some non-compact cases are given.  
                              
In the case of differentiable definable functions in an o-minimal structure and $U$ is bounded set (see \cite{Kur}), the author proved the \L ojasiewicz gradient inequality and the authors in \cite{BDL} proved it in the case of subanalytic functions. With some specific cases of o-minimal structures, other \L ojasiewicz-type inequalities was given in \cite{Loi}.

On the other hand, the classical \L ojasiewicz inequality has the relation with error bounds in Optimization. Let $f : \mathbb{R}^n \to \mathbb{R}$ be a continuous real-valued function. Set
\begin{equation}\label{Eqn1}
S:=\{x \in \mathbb{R}^n | f(x) \le 0\},
\end{equation}
and set $[f(x)]_+ := \max\{0, f(x)\}$.

We say that (\ref{Eqn1}) has a global H\"{o}lderian error bound if there exist $c >0, \alpha >0, \beta >0$ such that 
\begin{equation}\label{Eqn2}
d(x,S) \le c([f(x)]_+^\alpha + [f(x)]_+^\beta)
\end{equation}
for all $x \in \mathbb{R}^n$, where $d(x,S)$ denotes the Euclidean distance between $x$ and $S$. If, in addition, that $\alpha = \beta = 1$, then we refer (\ref{Eqn2}) as a global Lipschitzian error bound.

Note that $[f(x)]_+ = 0$ if and only if $x \in S$. Hence the existence of the \L ojasiewicz inequality with $[f(x)]_+$ over $K=\mathbb{R}^n$ is equivalent to the existence of the global H\"{o}lderian error bound of $S$.

In the convex case, the first results of error bounds was obtained in the work of many authors \cite{H}, \cite{R}, \cite{M}, \cite{AC}, \cite{KL}, \dots The existence of an error bound (Lipschitzian) usually requires the convexity and the so-called Slater condition. When the Slater condition is not satisfied and the set $S$ is defined by one or many polynomial inequalities, global H\"{o}lderian error bounds have been shown in \cite{LiG}, \cite{LL}, \cite{LS}, \cite{Y}, \dots 

In the non-convex case, the global H\"{o}lderian error bound for polynomial of degree $2$ was given in \cite[Theorem 3.1]{LS}. As far as we know, this is the first result, where a global H\"{o}lderian error bound for a \textit{non-convex} polynomial was established.

Recently, the author in \cite[Theorem A]{Ha} gave a criterion for the existence of a global H\"{o}lderian error bound (\ref{Eqn2}) in the case of polynomial of any degree, without the assumption the convexity and the Slater condition. Moreover, the author pointed out that if a polynomial satisfies the Palais-Smale condition then there exists a global H\"{o}lderian error bound.

In this paper, we will give some \L ojasiewicz-type inequalities. We will extend some results of \cite{Ha} from polynomial functions to continuous definable functions in an o-minimal structure. We also do not require functions to either be convex or satisfy the Slater condition. On the other hand, we will establish the \L ojasiewicz gradient inequality in a non-compact case with differentiable definable real-valued functions in an o-minimal structure.  

The rest of the paper is organized as follows. In Section 2, we recall a short introduction to o-minimal structures and some their properties. In Section 3, a criterion for the existence of \L ojasiewicz-type inequalities and \L ojasiewicz inequality of gradient will be proved. In Section 4, we give a necessary and sufficicent condition for which a global H\"{o}lderian error bound exists; moreover, a relation between the Palais-Smale condition and the existence of error bounds will be established in the end. 
\section{Preliminaries}
In this section, we recall some notions and results of geometry of o-minimal structures, which can be found in \cite{DM, D, C}.

\begin{definition}{\rm A structure expanding the real field $(\mathbb{R}, +, .)$ is a collection $\mathcal{O} = (\mathcal{O}_n)_{n \in \mathbb{N}}$ where each $\mathcal{O}_n$ is a set of subsets of the affine space $\mathbb{R}^n$, satisfying the following axioms:
\begin{enumerate}
  \item[1.] All algebraic subsets of $\mathbb{R}^n$ are in $\mathcal{O}_n$.
  \item[2.] For every $n$, $\mathcal{O}_n$ is closed under finite set-theoretical operations.
  \item[3.] If $A \in \mathcal{O}_n$ and $B \in \mathcal{O}_m$, then $A \times B \in \mathcal{O}_{m+n}$.
  \item[4.] If $\pi : \mathbb{R}^{n+1} \to \mathbb{R}^n$ is the projection on the first $n$ coordinates and $A \in \mathcal{O}_{n+1}$ then $\pi(A) \in \mathcal{O}_n$.\\
  The elements of $\mathcal{O}_n$ are called the {\em definable subsets} of $\mathbb{R}^n$. Moreover, if $\mathcal{O}$ satisfies:
  \item[5.] The elements of $\mathcal{O}_1$ are precisely the finite unions of points and intervals.\\
Then $\mathcal{O}$ is called an {\em o-minimal structure} on $\mathbb{R}$.
\end{enumerate}
}\end{definition}

\begin{example}{\rm
A semi-algebraic set is finite union of sets $S = \{x \in \mathbb{R}^n | f(x) = 0, g_j(x) < 0, j = 1, \dots, m\}$ where $f, g_j$ are polynomials in $\mathbb{R}[x_1, \dots, x_n]$.\\
The collection $\mathcal{O}$ of all semi-algebraic sets in $\mathbb{R}^n$ for all $n \in \mathbb{N}$ is an o-minimal structure on $\mathbb{R}$.    
}\end{example}
Perhaps the writing down projections in order to show that a subset is definable will be boring. We are more used to write down formulas. Let us specify what is meant \textit{first-order formula} (of the language of the o-minimal structure). A first-order formula is constructed according to the following rules.
\begin{enumerate}
\item If $P \in \mathbb{R}[X_1, \dots, X_n]$, then $P(X_1, \dots, X_n) = 0$ and $P(X_1, \dots, X_n) > 0$ are first-order formulas.
\item If $A$ is a definable subset of $\mathbb{R}^n$, then $x \in A$ (where $x = (x_1, \dots, x_n)$) is a first-order formula.
\item If $\Phi(x_1, \dots,x_n)$ and $\Psi(x_1, \dots, x_n)$ are the first-order formulas, then \{$\Phi$ and $\Psi$\}, \{$\Phi$ or $\Psi$\}, \{not $\Phi$\}, \{$\Phi \Rightarrow \Psi$\} are first-order formulas.
\item If $\Phi(y, x)$ is a first-order formula (where $y = (y_1,\dots, y_p)$ and $x = (x_1, \dots, x_n)$) and $A$ is a definable subset of $\mathbb{R}^n$, then $\exists x \in A\ \Phi(y, x)$ and $\forall x \in A\ \Phi(y, x)$ are first-order formulas.
\end{enumerate}
\begin{theorem}[\cite{C}, Theorem 1.13]\label{first-order} If $\Phi(x_1, \dots, x_n)$ is a first-order formula, the set of $(x_1, \dots, x_n)$ in $\mathbb{R}^n$ which satisfy $\Phi(x_1, \dots,x_n)$, is definable.
\end{theorem}
\begin{remark}{\rm By the rule (4) and the above theorem, the sets $\{x \in \mathbb{R}^n: \exists x_{n+1} (x, x_{n+1})\in A\}$ (image of $A$ by projection) and $\{x \in \mathbb{R}^n: \forall x_{n+1} (x, x_{n+1}) \in A\}$ (complement of the image of the complement of $A$ by projection) are definable.
}\end{remark}
\begin{definition}{\rm A map $f: A \to \mathbb{R}^p$ (where $A \subset \mathbb{R}^n$) is called {\em definable} if its graph is a definable subset of $\mathbb{R}^n \times \mathbb{R}^p$.}
\end{definition}
With any o-minimal structure, we have some elementary properties
\begin{proposition}\label{Prop0}\begin{enumerate}
\item[(i)] The closure, the interior and the boundary of a definable set are definable.
\item[(ii)] Compositions of definable maps are definable.
\item[(iii)] Images and inverse images of definable sets under definable maps are definable.
\item[(iv)] Infimum of a bounded below definable function and supremum of a bounded above definable function are definable functions.
\end{enumerate}
\end{proposition}
The reader can be found the proofs of these properties in \cite{DM, D}.
\begin{proposition}
If function $f: \mathbb{R}^n \to \mathbb{R}$ is definable then the set $S = \{x \in \mathbb{R}^n | f(x) \le 0\}$ is definable.
\end{proposition}
\begin{proof}
By definition, $\Gamma_f = \mathbb{R}^n \times f(\mathbb{R}^n)$ is definable.\\
Let consider the following projection \begin{align*}
\pi: \mathbb{R}^{n+1} &\to \mathbb{R},\\ 
(x_1, \dots, x_n, x_{n+1}) &\mapsto x_{n+1}. 
\end{align*}
By the definition of first-order formula, the set $\pi(\Gamma_f) = \{y \in \mathbb{R} | y =f(x),\ \text{for some}\ x \in \mathbb{R}^n\}$ is definable. Similarly, the set $\{y \in \mathbb{R} | y \le 0\}$ is definable.\\
So $S = \pi(\Gamma_f) \cap \{y \le 0\}$ is definable.
\end{proof}

\begin{proposition}\label{prop1}
If $S$ is a definable set and $S \ne \emptyset$ then the function $d: \mathbb{R}^n \to \mathbb{R}$ defined by $$d(x,S) = \inf_{y \in S}\|x- y\|$$ is well-defined and is a definable function; moreover, it is a continuous function on $\mathbb{R}^n$.
\end{proposition}
\begin{proof}
The set $\{\|x - y\| : y \in S\}$ is an image of $S$ by the definable function $y \mapsto \|x - y\|$, so it is definable subset. Since $S \ne \emptyset$, $d$ is well-defined.\\
Let consider its graph, $\Gamma_d = \{(x,t) \in \mathbb{R}^{n+1} | t \ge 0\ \text{and}\ \forall y \in S: t^2 \le \|x - y\|^2\ \text{and}\ \forall \epsilon \in \mathbb{R}, \epsilon > 0 \Rightarrow \exists y \in S: t^2 + \epsilon > \|x - y\|^2\}$.\\
This set is definable because it is defined by first-order formulas. Hence $d(x,S)$ is a definable function.\\
By the triangle inequality, we have $|d(x, S) - d(x_0, S)| \le d(x, x_0)$. Therefore $x \to x_0$ implies $d(x, S) \to d(x_0, S)$. Hence $d(x, S)$ is a continuous function.
\end{proof}
\begin{proposition}\label{prop2}
Let $f : \mathbb{R}^n \to \mathbb{R}$ be a differentiable, definable function in some o-minimal structure. Then $\partial f/ \partial x_j, j =1, \dots, n$ are definable functions and $\nabla f(x)$ (gradient of $f$) is an definable mapping. 
\end{proposition}
\begin{proof}
By the definition of partial derivatives, we have $\partial f/\partial x_j$ are defined by $$\partial f/\partial x_j (a) = \lim\limits_{x_j \to a_j}\dfrac{f(x_1, \dots, x_j, \dots, x_n) - f(a_1, \dots, a_j, \dots, a_n)}{x_j - a_j}, a\in\mathbb{R}^n,$$ so we have
$$- \epsilon < \dfrac{f(x_1, \dots, x_j + h, \dots, x_n) - f(x_1, \dots, x_n)}{h} - \partial f/\partial x_j < \epsilon, \forall \epsilon > 0, h > 0, j = 1, \dots, n.$$ This is a first-order formula. By Theorem \ref{first-order}, $\partial f/\partial x_j$ is definable function. This implies that $\nabla f(x)$ is definable. 
\end{proof}
\noindent The following useful result is a property of semialgebraic functions in one variable.
\begin{lemma}[\cite{DM}, growth dichotomy Lemma] \label{GrowthDichotomyLemma}
Let $f \colon (0, \epsilon) \rightarrow {\Bbb R}$ be a semi-algebraic function with $f(s) \ne 0$ for all $s \in (0, \epsilon).$ Then 
there exist constants $c \ne 0$ and $q \in {\Bbb Q}$ such that $f(s) = cs^q + o(s^q)$ as $s \to 0^+.$
\end{lemma}
 
\noindent The following property is important to our purpose.
\begin{theorem}[monotonicity theorem] Let $f : (a,b) \to \mathbb{R}$ is a definable function, $-\infty \le a < b \le +\infty$. Then there exist $a_0, a_1, \dots, a_{k+1}$ with $a = a_0 < a_1 < \dots < a_k < a_{k+1} = b$ such that $f$ is continuous on each interval $(a_i, a_{i+1})$, moreover $f$ is either strictly monotone or constant on each $(a_i, a_{i+1}), i = 1, \dots, k$.
\end{theorem}
\noindent The proof of this theorem can be found in \cite{DM, D, C}.

We now recall notion of the subdifferential of a continuous function. This notion plays the role of the usual gradient map, which can be found in \cite{RW, Cl}.
\begin{definition}{\rm
\begin{enumerate}
  \item[(i)] The {\em Fr\'echet subdifferential} $\hat{\partial} f(x)$ of a continuous function $f \colon {\Bbb R}^n \rightarrow {\Bbb R}$ at $x \in {\Bbb R}^n$ is given by $$\hat{\partial} f(x) := \left \{ v \in {\Bbb R}^n \ | \ \liminf_{\| h \| \to 0, \ h \ne 0} \frac{f(x + h) - f(x) - \langle v, h \rangle}{\| h \|} \ge 0 \right \}$$
  \item[(ii)]  The {\em limiting subdifferential} at $x \in {\Bbb R}^n,$ denoted by ${\partial} f(x),$ is the set of all cluster points of sequences $\{v^k\}_{k \ge 1}$ such that $v^k\in \hat{\partial} f(x^k)$ and $(x^k, f(x^k)) \to (x, f(x))$ as $k \to \infty.$
\end{enumerate}
}\end{definition}

\begin{remark}\label{remark1}{\rm \begin{description}
\item[(i)] It is easy to show that for a continuous function $f$ on $\mathbb{R}^n$, the set $\{x: \hat{\partial} f(x) \ne \emptyset\}$ is dense set in $\mathbb{R}^n$.
\item[(ii)] It is not hard to show that if $f$ is a definable function then $\hat{\partial} f(x)$ and $\partial f(x)$ are definable sets (\cite[Prop 3.1]{I1}).
\end{description} 
}\end{remark}

\begin{definition}{\rm
By using the limiting subdifferential $\partial f,$ we define the {\em
nonsmooth slope} of $f$ by
$${\frak m}_f(x) := \inf \{ \|v\| : v \in {\partial} f(x) \}.$$
By definition, ${\frak m}_f(x) = + \infty$ whenever ${\partial} f(x)
= \emptyset.$
}\end{definition}

\begin{definition}{\rm
The \textit{strong nonsmooth slope} of function $f$ is defined as follows $$|\nabla f|(x) : = \lim_{h \to 0}\sup_{h \ne 0}\frac{[f(x) - f(x + h)]_+}{\|h\|},$$ with $[a]_+ = \max\{a, 0\}$.
}\end{definition}

The relationship between nonsmooth slope, strong nonsmooth slope and subdifferential is following (see for details in \cite{I2}): $$ \inf\{\|y\|: y \in \hat{\partial}f(x)\} \ge |\nabla f|(x) \ge \mathfrak{m}_f(x). $$

\begin{remark}\label{remark2}{\rm \begin{description}
\item[(i)] It is not hard to show that if $f$ is a definable function then $\mathfrak{m}_f(x)$ and $|\nabla f|(x)$ are definable (\cite[Prop 3.1]{I1}).
\item[(ii)] If $f$ is a differentiable function then the above notions coincide with the usual concept of gradient; that is: $\partial f(x) = \hat{\partial}f(x) = \{\nabla f(x)\}$ and hence $\mathfrak{m}_f(x) = |\nabla f|(x) = \|\nabla f(x)\|$.
\end{description}
}\end{remark}
\section{Main results}
\subsection{\L ojasiewicz-type inequalities}\quad\\
The following results extend the results of \cite{Ha} (see also \cite{DHN}) from polynomial functions to continuous definable functions. The proof follows the steps of proofs of Theorem 2.1 and 2.2 in \cite{Ha}, but we use monotonicity theorem instead of growth dichotomy lemma.

\begin{proposition}[\L ojasiewicz-type inequality "near to the set $S$"]\label{near}
Let $f \colon {\Bbb R}^n \rightarrow {\Bbb R}$ be a continuous definable function. Assume that $S := \{x \in {\Bbb R}^n \ | \ f(x) \le 0 \} \ne \emptyset.$ Let $[f(x)]_+ := \max\{f(x), 0\}.$ Then the following two statements are equivalent.
\begin{enumerate}
\item[{\rm(i)}] For any sequence $x^k \in \R^n \setminus S,$ with $x^k \to \infty,$ it holds that 
$$f(x^k) \to 0 \quad \Longrightarrow \quad d(x^k,S) \to 0;$$
\item[{\rm(ii)}] There exist $\delta > 0$ and a function $\mu: [0,\delta] \to \mathbb{R}$ which is definable, continuous and strictly increasing on $[0,\delta)$ with $\mu(0) = 0$ such that $$ \mu([f(x)]_+) \ge d(x,S),\quad \forall x \in f^{-1}((-\infty, \delta]). $$
\end{enumerate}
\end{proposition}

\begin{proof}\quad \\ 
\noindent$(ii) \Rightarrow (i):$ Assume that $x^k \not\in S,$ $x^k \to \infty$ and $f(x^k) \to 0.$ 
We have $[f(x^k)]_+ = f(x^k)$. By the continuity of $\mu$ at $0$, we get $\mu(f(x^k)) \to 0$. Note that $0 < f(x^k) < \delta$ if $k \gg 1.$ Then it follows from the inequality in $(ii)$ that $d(x^k, S) \to 0$.\\

\noindent $(i) \Rightarrow (ii):$ Without loss of generality, we can suppose that $S \ne \mathbb{R}^n$. Then there exists $t_0 > 0$ such that $f^{-1}(t_0) \ne \emptyset.$ Because $f$ is continuous, $f^{-1}(t) \ne \emptyset$ for all $0 \le t \ll 1.$ 

Let $\mu(t) := \sup\limits_{x \in f^{-1}(t)}d(x, S), t \ge 0.$ We will show that there exists $\delta > 0$ sufficient small such that $\mu(t)$ have desired properties. Clearly, $\mu(0) = 0$.

We now show that there exists $\delta > 0$ such that $\mu(t) < +\infty$ for all $t \in [0, \delta)$. By contradiction, assume that there exists a sequence $t_k > 0, t_k \to 0,$ such that $\mu(t_k)  = \infty$ for all $k.$ This implies the existence of sequence $x^k \in f^{-1}(t_k)$ such that $d(x^k, S) \to +\infty$ as $k \to \infty.$ Hence $x^k \to \infty.$ Contradiction.

So $\mu(t) < +\infty$ on $[0, \delta]$ with $\delta > 0$. By Proposition \ref{prop1} and Proposition \ref{prop2}((iv)), we have $\mu(t)$ is definable on $[0, \delta].$

Using the monotonicity theorem, the function $\mu$ is continuous and monotone on $(0, \delta]$ if $0 < \delta \ll 1.$ 

We now show that $\mu$ is continuous at $0.$ Suppose $\mu$ is not continuous at $0$. That means there exists a sequence $t_k \to 0$ such that $\mu(t_k) = \sup\limits_{x \in f^{-1}(t_k)}d(x, S) \nrightarrow 0$. Hence, there exists a sequence $x^k \in f^{-1}(t_k)$ such that $t_k = f(x^{k}) \to 0$ and $d(x^{k}, S) \nrightarrow 0.$ On the other hand, $x^k \to \infty.$ Indeed, if there exists $x < \infty$ such that $x^k \to x$ then by the continuity of $f$, $f(x^{k}) \to f(x)$, this implies $f(x) = 0$. That means $d(x^k, S) \to 0$, contradiction. So we have a sequence $x^k \to \infty, f(x^k) \to 0$ and $d(x^k, S) \nrightarrow 0$. This contradicts $(i)$.

Hence $\mu$ is continuous and monotone on $[0, \delta].$

Note that by $\mu(0) = 0$ and $\mu(t) > 0, \forall t \in (0, \delta)$, if $\delta$ is sufficient small then $\mu(t)$ is strictly increasing on $[0, \delta].$

For $0 < t < \delta$, let $x \in f^{-1}(t)$, then we have $\mu(t) = \sup\limits_{a \in f^{-1}(t)}d(a, S) \ge d(x, S)$.
          
Hence $\mu([f(x)]_+) \ge d(x, S), \forall x \in f^{-1}((-\infty, \delta])$.
\end{proof}

\begin{remark}{\rm Note that the condition that $\mu$ is continuous at $0$ and $\mu(0) = 0$ in $(ii)$ is necessary.\\
Let us consider the function $f: \mathbb{R} \to \mathbb{R}, x \mapsto \dfrac{x}{1+x^2}$. The function $f$ is a differentiable semialgebraic function because its graph is the set $\{(x, y) \in \mathbb{R}^2|  (1+x^2)y = x\}$. Then $f$ is a definable function.\\
We have $S = (-\infty, 0]$. Then we choose $\mu(t) :=\sup\limits_{\frac{x}{1+x^2} = t}d(x, S)$ on $0<t<\dfrac{1}{2}$. This function is definable, continuous on $(0, \dfrac{1}{2})$ but not continuous at $0$.\\
Moreover, $x^k \to +\infty$ satisfies $f(x^k) \to 0$ but $d(x^k, S) \to +\infty$, so the statement $(i)$ fails.
}\end{remark}
\begin{proposition}[\L ojasiewicz-type inequality "far from the set $S$"]\label{far} 
Suppose that for any sequence $x^k \in \R^n \setminus S,$ with $x^k \to \infty$ and
$$d(x^k,S) \to \infty \ \text{we have}\ f(x^k) \to \infty;$$
Then there exist $r > 0$ and a function $\mu: [r,+\infty) \to \mathbb{R}$ which is definable, increasing and continuous on $[r,+\infty)$ such that 
$$ \mu([f(x)]_+) \ge d(x,S), \forall x \in f^{-1}([r,+\infty)). $$
\end{proposition}
\begin{proof}\quad \\
Let us consider two cases:\\
\indent Case 1. The function $f$ is bounded from above, i.e. $r := \sup_{x \in {\Bbb R}^n} f(x) < +\infty.$ 

By the assumption, there exists $M > 0$ such that $d(x, S) \le M$ for all $x \in \mathbb{R}^n$. For all $x \in f^{-1}([r', r))$ ($0< r' < r$), $$f(x) \ge r' = \frac{r'}{M} M \ge \frac{r'}{M}d(x,S),$$ 
Then the function $\mu(t) := \frac{M}{r'} t$ with $t \ge r'$ have required properties.

Case 2. The function $f$ is not bounded from above. By continuity of $f$ and $S \ne \emptyset,$ we have $f^{-1}(t) \ne \emptyset$ for all $t \ge 0.$ Set $\mu(t) = \sup\limits_{x \in f^{-1}(t)}d(x, S)$.\\
We claim that there exists $r \gg  1$ such that $\mu(t) = \sup\limits_{x \in f^{-1}(t)}d(x, S) < \infty$ for all $t \ge r.$ By contradiction, assume that
$\mu(t) = \infty$ for some $t \gg 1$. Then there exists a sequence $x^k \in f^{-1}(t)$ such that $d(x^k, S) \to \infty.$ Of course $x^k \to \infty$, this contradicts the assumption.

So $\mu(t) < +\infty, \forall t \in [r, +\infty).$ This implies that $\mu$ is a definable function on $[r, +\infty)$. By monotonicity theorem, $\mu$ is continuous and monotone on $[r, +\infty)$ for $r \gg 1$.

Let
$$M := \sup\limits_{t \in [r, +\infty)} \mu(t).$$
We have two subcases:

Case 2.1. $M = +\infty.$ Then $\lim_{t \to +\infty} \mu(t) = +\infty.$ This means that for $r \gg 1,$ the function $\mu$ is strictly increasing on $[r, +\infty).$ Furthermore
$$\mu([f(x)]_+) =  \mu(f(x)) \ge d(x, S), \quad \forall x \in f^{-1}([r,+\infty)).$$

Case 2.2. $M < +\infty.$ Then, for all $x$ such that $f(x) \ge r$ we have $d(x, S) \le M$, therefore
$$f(x) \ge r = \frac{r}{M} M \ge \frac{r}{M} d(x, S).$$
The function $\mu := \frac{M}{r}t, t \ge r,$ has required properties.
\end{proof}

\begin{remark}{\rm Note that the converse of the above theorem is false.\\
Indeed, consider the function $f: \mathbb{R} \to \mathbb{R}, x \mapsto \dfrac{x}{\sqrt{1+x^2}}$. The function $f$ is a differentiable semialgebraic function since its graph is the set $\{(x, y) \in \mathbb{R}^2|  (1+x^2)y^2 = x^2\} \cap \{xy > 0\}$.\\
We have $S = (-\infty, 0]$. We choose $0<r <1$ and let $\mu(t) :=\begin{cases}\sup\limits_{\frac{x}{\sqrt{1+x^2}} = t}d(x, S) & \text{on } [r, 1)\\ +\infty & \text{on } [1, +\infty) \end{cases}$. This function is definable, increasing and continuous.\\
In the other hand, we have $x^k \to +\infty$, $d(x^k, S) \to +\infty$ and $f(x^k) \to 1$.\\
}\end{remark}
\subsection{Global H\"{o}lderian error bound for continuous definable functions in o-minimal structures}\quad \\
The following criterion extends the error bound result of \cite{Ha} from polynomial functions to definable functions in o-minimal structures.
\begin{theorem}\label{Global}
Let $f \colon {\Bbb R}^n \rightarrow {\Bbb R}$ be a continuous definable function. Assume that $S := \{x \in {\Bbb R}^n \ | \ f(x) \le 0 \} \ne \emptyset$ and $[f(x)]_+ := \max\{f(x), 0\}.$ Then the following two statements are equivalent
\begin{enumerate}
\item[{\rm (i)}] For any sequence $x^k \in \mathbb{R}^n \setminus S, x^k \to \infty$, we have
\begin{enumerate}
\item[(i1)] if $f(x^k) \to 0$ then $d(x^k,S) \to 0;$
\item[(i2)] if $d(x^k,S) \to \infty$ then $f(x^k) \to \infty.$
\end{enumerate}

\item[{\rm (ii)}] There exists a function $\mu: [0,+\infty) \to \mathbb{R}$, which is definable, strictly increasing and continuous on $[0, +\infty)$ with $\mu(0) = 0, \lim\limits_{t \to +\infty} \mu(t) = +\infty,$ such that 
$$ d(x,S) \le \mu([f(x)]_+),\quad \forall x \in \mathbb{R}^n. $$
\end{enumerate}
\end{theorem}
\begin{proof}
The implication (ii) $\Rightarrow$ (i) is straightforward. We prove the implication (i) $\Rightarrow$ (ii). 

Indeed, by Theorems \ref{near} and \ref{far}, there exist two continuous, strictly increasing, definable functions $\mu_1$ on $[0, \delta]$ with $0 < \delta \ll 1$ and $\mu_2$ on $[r, +\infty)$ with $r \gg 1$ such that

$$ d(x,S) \le \mu_1([f(x)]_+),\forall x \in f^{-1}((-\infty, \delta]). $$
and
$$d(x,S) \le \mu_2([f(x)]_+),\forall x \in f^{-1}([r,+\infty)).$$

On the other hand, by assumption $(i2)$, there exists $M > 0$ such that $d(x, S) \le M$ for all $x \in f^{-1}([\delta, r]).$ Then
$$f(x) \ge \delta = \frac{\delta}{M} M \ge \frac{\delta}{M} d(x, S)$$ 
for all $x \in f^{-1}([\delta, r]).$ Put $\mu_3(t) := \frac{M}{\delta} t$ with $t \in [\delta, r],$ we get $\mu_3(t) \ge d(x, S)$ and $\mu_3$ is a increasing function on $[\delta, r]$.

By definition of $\mu_3$ and $\lim\limits_{t \to 0}\mu_1(t) = 0$ (Theorem \ref{near}), we may choose $\delta$ such that $\mu_1(t) \le \mu_3(\delta) = M, \forall t \in [0, \delta]$. Indeed, if $\exists t \in [0, \delta]$ such that $\mu_1(t) > \mu_3(\delta)$, then we put $$ M' := \max\{\sup\limits_{t \in [0, \delta]}\mu_1(t), M\} \ \text{and}\  \mu_3(t) := \frac{M'}{\delta}t, $$ so we have $\mu_1(t) \le M' = \mu_3(\delta), \forall t \in [0, \delta]$.

Similarly, by definition of $\mu_3(t)$ and $\mu_2(t)$, we may choose $r$ such that $\mu_3(r) = \frac{M}{\delta}r \le \mu_2(t), \forall t \in [r, +\infty)$. Indeed, if $\exists t \in [r, +\infty)$ such that $\mu_3(r) > \mu_2(t)$, then we may choose $\mu_2'(t) \ge \mu_2(t) + C$ with $C = \mu_3(r)$, so we have $d(x, S) \le \mu_2(t) < \mu_2'(t), t \in [r, +\infty)$ and $\mu_3(r) \le \mu_2'(t), \forall t \in [r, +\infty)$. Moreover, by definition of $\mu_2$, we may choose $\mu_2'(t)$ as above such that if $r \gg 1$ then $\mu_2'(t)$ is strictly increasing on $[r, +\infty)$ and $\lim\limits_{t \to +\infty}\mu_2'(t) = +\infty$.

Combine three functions $\mu_1, \mu_2, \mu_3$ and note that we may choose suitable $\delta, r$ and $M$ as above, we get the function 
$\mu(t) = \begin{cases}\mu_1(t) & \forall t \in [0, \delta]\\
\mu_3(t) & \forall t \in [\delta, r]\\
\mu_2'(t) & \forall t \in [r, +\infty) \end{cases}.$ The function $\mu$ is definable, strictly increasing and continuous and $\mu$ satisfies $(ii)$.
\end{proof}
\subsection{The relation between the Palais-Smale condition and the existence of error bounds}\quad\\
In this section, we consider continuous functions in an o-minimal structure.
\begin{definition}{\rm 
Given a continuous function $f \colon {\Bbb R}^n \rightarrow {\Bbb R}$ and a real number $t,$ we say that $f$ satisfies the {\em Palais-Smale condition at the level} $t,$ if every sequence $\{x^k\}_{k \in {\Bbb N}} \subset {\Bbb R}^n$ such that $f(x^k) \to t$ and ${\frak m}_f(x^k) \to 0$ as $k \to \infty$ possesses a convergence subsequence.
}\end{definition}

The following theorem also extends Theorem B in \cite{Ha} from polynomial functions to continuous definable functions.

\begin{theorem} \label{Palais-Smale}
Let $f \colon {\Bbb R}^n \rightarrow {\Bbb R}$ be a continuous definable function. Assume that $S := \{x \in {\Bbb R}^n \ | \ f(x) \le 0 \} \ne \emptyset.$ If $f$ satisfies the Palais-Smale condition at each level $t \ge 0$, then there exists a function $\mu: [0,+\infty) \to \mathbb{R}$, which is definable, strictly increasing and continuous $\mu(0) = 0, \lim_{t \to \infty} \mu(t) = \infty,$ such that 
$$ d(x,S) \le \mu([f(x)]_+),\quad \forall x \in \mathbb{R}^n. $$
\end{theorem}

\begin{proof}
By Theorem \ref{Global}, it is enough to show that $f$ satisfies the Palais-Smale condition at each value $t \ge 0$, then there is no sequence $x^k \to \infty, x^k \in \mathbb{R}^n \setminus S$ such that $$ f(x^k) \to 0 \ \text{but}\ d(x,S) > \delta > 0$$ or $$ d(x^k, S) \to \infty \ \text{but}\ 0 \le f(x^k) \le M $$ for some $\delta > 0$ and $M > 0$. In case of continuous definable functions, we use the subdifferential instead of the gradient in \cite{Ha}.

By contradiction, first of all, assume that for a sequence $x^k \to \infty, x^k \in \mathbb{R}^n \setminus S$, we have $f(x^k) \to 0$ and $d(x^k , S) \ge \delta > 0$. Similarly to the proof of \cite[Theorem B]{Ha}, by using Ekeland Variational Principle (\cite{E}), we obtain a sequence $y^k$ such that $$ \frac{1}{\|h\|}(f(y^k + h) - f(y^k))\ge -\sqrt{\epsilon_k} $$ with $h \in \mathbb{R}^n, 0 < \|h\| < \dfrac{\delta}{2}$ and $\epsilon_k= f(x^k)$. This implies that
$$ \frac{1}{\|h\|}( f(y^k) - f(y^k + h)) \le \sqrt{\epsilon_k},$$
or
$$ \frac{1}{\|h\|} [f(y^k) - f(y^k + h)]_+ \le \sqrt{\epsilon_k}.$$
By the definition of the strong slope, we have 
$$0 \le  |\nabla f|(y^k)  = \limsup_{h \to 0, h \ne 0} \dfrac{[f(y^k) - f(y^k + h)]_+}{\|h\|} \le \sqrt{\epsilon_k}.$$ Thus 
$$ 0 \le \mathfrak{m}_f(y^k) \le |\nabla f|(y^k) \le \sqrt{\epsilon_k}.$$
Letting $k \to \infty$, we get $\mathfrak{m}_f(y^k) \to 0$. So we have found a sequence $y^k \to \infty, y^k \in \mathbb{R}^n \setminus S, \mathfrak{m}_f(y^k) \to 0$ and $f(y^k) \to 0$. This means that $f$ does not satisfy the Palais-Smale condition at the value $t=0$, a contradiction. So we get (i1) of Theorem \ref{Global}.

Now, suppose that for some sequence $x^k \in \mathbb{R}^n \setminus S$ with $x^k \to \infty$ such that 
$$ d(x^k, S) \to \infty\ \text{and}\ f(x^k) \nrightarrow \infty. $$
Without loss of generality, we may assume that $f(x^k) \to t_0$ with $t_0 \in [0, + \infty)$. Again, by the similar arguments as in \cite[Theorem B]{Ha}, we have a sequence $y^k$ such that $0 < f(y^k) \le f(x^k)$ and $$ \frac{1}{\|h\|}(f(y^k + h) - f(y^k))\ge -\epsilon_k\cdot \lambda_k $$ with $h \in \mathbb{R}^n, 0 < \|h\| < \dfrac{\delta}{2}, \epsilon_k= f(x^k)$ and $\lambda_k = \dfrac{2}{d(x^k,S)}$. This implies that $$ \frac{1}{\|h\|} [f(y^k) - f(y^k + h)]_+ \le \epsilon_k\lambda_k. $$
By the definition of the strong slope, we have $$ 0 \le \mathfrak{m}_f(y^k) \le |\nabla f|(y^k) \le \epsilon_k\lambda_k = \dfrac{2\epsilon_k}{d(x^k, S)}.$$
Letting $k \to \infty$ we have $\epsilon_k = f(x^k) \to t_0$ and $d(x^k, S) \to \infty$. Therefore $\mathfrak{m}_f(y^k) \to 0$.\\
Consequently, since $0 < f(y^k) \le f(x^k)$, $y^k$ has a subsequence $y'^k$ such that $f(y'^k) \to t_1$ with $0 \le t_1 \le t_0$ which satisfies $$ y'^k \to \infty, \mathfrak{m}_f(y'^k) \to 0\ \text{and}\ f(y'^k) \to t_1. $$ This means that $f$ does not satisfy the Palais-Smale condition at $t_1$, contradiction. So we get (i2) of Theorem \ref{Global}. The theorem is proved.
\end{proof}

\subsection{A nonsmooth slope inequality near the fiber for continuous definable functions in an o-minimal structure}\quad\\
In the case $U$ is not bounded set, the classical \L ojasiewicz gradient inequality is not always true. We can see it in the following example
\begin{example}\label{Ex2}{\rm Consider the following example: $ f(x,y) = (xy - 1)^2 + (x - 1)^2$ and $U = \mathbb{R}^2$. Let $x^k$ be $(\dfrac{1+ k}{1 + k^2}, k)$, we have:
\begin{itemize}
\item $x^k \to \infty$.
\item $\nabla f(x^k) = (0, 2\dfrac{k^2 - 1}{(1 + k^2)^2}) \to 0$.\\
But
\item $f(x^k) = (\dfrac{(1 + k).k}{1 + k^2} - 1)^2 + (\dfrac{1 + k}{1 + k^2} - 1)^2 \to 1.$
\end{itemize}
We prove that $\nexists \delta > 0, C > 0, \rho \in \mathbb{R}$ such that $\|\nabla f(x)\| \ge C|f(x)|^\rho\ \text{for}\ x \in f^{-1}(D_\delta)$ with $D_\delta = \{t:|t| < \delta\}$ By contradiction, assume that there are $\delta > 0$, $ C > 0$ and $\rho \in\mathbb{R}$ such that the \L ojasiewicz gradient inequality holds.  We see that $\nabla f(\dfrac{1}{k}+1, \dfrac{k}{k+1}) \to 0$ and $f(\dfrac{1}{k}+1, \dfrac{k}{k+1}) \to 0$. Hence $\rho > 0$. On the other hand, $\nabla f(\dfrac{1+ k}{1 + k^2}, k) \to 0$ and $f(\dfrac{1+ k}{1 + k^2}, k) \to 1$; so $\rho \le 0$, contradiction.
}\end{example}
We shall give a criterion for the existence of \L ojasiewicz nonsmooth slope inequality on $f^{-1}(D_\delta)$.

Let $\widetilde{K}(f) := \{ t \in \mathbb{R} |\ \exists x^k, {\frak m}_f(x^k) \to 0, f(x^k) \to t\}$ and we call it the set of asymptotic critical values.
\begin{theorem}\label{gradient} Let $f :\mathbb{R}^n \to \mathbb{R}$ be a continuous definable function in some o-minimal structure and suppose that $\widetilde{K}(f) \cap D_\delta = \{0\}$. Then the following two statements are equivalent.
\begin{enumerate}
\item[(i)] For any sequence $x^k \to \infty$, $\frak{m}_f(x^k) \to 0$ implies $f(x^k) \to 0$.
\item[(ii)] There exists a function $\varphi: (0,\delta) \to \mathbb{R}$, which is definable, monotone and continuous such that 
$$ \frak{m}_f(x) \ge \varphi(|f(x)|),\quad \forall x \in f^{-1}(D_\delta). $$  
\end{enumerate} 

\end{theorem}

\begin{proof}[Proof of Theorem \ref{gradient}]\quad\\
$(i) \Rightarrow (ii):$ Let $\varphi(t) := \inf\{\frak{m}_f(x) : |f(x)| = t \}$, it is easy to see that $\varphi$ is a definable function (see Remark \ref{remark2} and Proposition \ref{Prop0}).

\noindent\textbf{Claim:} There exists $\delta_1$ such that $\varphi(t) > 0, \forall t \in (0,\delta_1)$.

Indeed, by the assumption $\widetilde{K}(f) \cap D_\delta = \{0\}$, $(0, \delta)$ has no critical point of $f$. Assume that there exists a value $t \in (0,\delta')$ such that $\varphi(t) = 0$. Then there exists a sequence $t_k$ such that $t_k \to t$ implies $\varphi(t_k) \to 0$. Therefore there exists a sequence $x^k$ such that $f(x^k) = t_k$ and $\frak{m}_f(x^k) \to 0$. So we have $\frak{m}_f(x^k) \to 0$ but $f(x_k) \to t \ne 0$, this contradicts with $(i)$. The claim is proved.

On the other hand, by monotonicity theorem, $\varphi(t)$ is continuous and monotone on $(0, \delta)$ for $0 < \delta \ll 1$.

By the definition of $\varphi$, we get $\varphi(t) \le \frak{m}_f(x), \forall x \in f^{-1}(D_\delta)$, which means that $\frak{m}_f(x) \ge \varphi(|f(x)|), \forall x \in f^{-1}(D_\delta).$  

\noindent$(ii) \Rightarrow (i):$ straightforward.
\end{proof}
\begin{remark}{\rm Cardinal of the set $\widetilde{K}_\infty(f)$ can be infinite. Indeed, consider the following example\\
Consider $f(x,y) = \dfrac{x}{1 +y^2}$ in the o-minimal structure of all semialgebraic sets, then any $t \in \mathbb{R}$ is belong to $\tilde{K}_\infty(f)$, by the sequence $x^k = (t(1+k^2), k)$. It is easy to see that $x^k \to \infty, \frak{m}_f(x^k) = \|\nabla f(x^k)\| = \sqrt{(\dfrac{1}{1+k^2})^2 + (\dfrac{2tk}{1+k^2})^2} \to 0$ and $f(x^k) = t$.
}\end{remark}

\begin{remark}{\rm In Theorem \ref{gradient}, if $f$ is a polynomial then $\varphi(t)$ is a semialgebraic function in one variable. By Growth Dichotomy Lemma, there exists $a > 0$ and $u \in \mathbb{R}, u > 0$ such that $$ \varphi(t) = at^u + o(t^u)\ \text{as} \ t \in (0, \epsilon), \epsilon \ll 1. $$ This implies $\varphi(t) \ge ct^u, \forall t \in (0, \epsilon)$. By definition of $\varphi$, we have $\|\nabla f(x)\| \ge \varphi(t) \ge ct^u$. Note that $t = |f(x)|$, so we get the \L ojasiewicz gradient inequality on $f^{-1}(D_\delta)$.
}\end{remark}
\subsubsection*{Acknowledgments}
The author would like to thank Prof. Ha Huy Vui for his proposing the problem, Assoc. Prof. Pham Tien Son for his suggestions and Dr. Dinh Si Tiep for his useful discussions.
\bibliographystyle{amsalpha}

\end{document}